\documentclass[12pt]{article}

\usepackage{fullpage}
\usepackage{epsfig}
\usepackage{color}
\usepackage{bbm}

\newtheorem{theorem}{Theorem}    % Specify Theorem
\newtheorem{lemma}{Lemma}      % Specify Lemma
\newtheorem{remark}{Remark}
\newtheorem{example}{Example}
\newtheorem{question}{Question}
\newtheorem{conjecture}{Conjecture}

\newtheorem{corollary}{Corollary}      % Specify corollary
\def\eod{\vrule height 6pt width 5pt depth 0pt}
\newenvironment{proof}{\noindent {\bf Proof:} \hspace{.2em}}
                      {\hspace*{\fill}{\eod}}
                      {}
\newcommand{\oa}{\overline{a}}
\newcommand{\ob}{\overline{b}}
\newcommand{\oc}{\overline{c}}
\newcommand{\PF}{\mathsf{PF}}
\newcommand{\comment}[1]{}
\newcommand{\inn}{\mathsf{invs}}
\newcommand{\Dist}{\mathsf{Dist}}
\newcommand{\IP}{\mathsf{IP}}
\newcommand{\NIP}{\mathsf{NIP}}
\newcommand{\RSn}{{\mathcal R}({\cal S}_n)}
\newcommand{\B}{{\mathcal B}}
\newcommand{\Sh}{{\mathcal S}}

\newcommand{\reals}{\mathbbm{R}}

\begin{document}

\title{Interpreting the two variable Distance enumerator of the 
Shi hyperplane arrangement}
\author{Sivaramakrishnan Sivasubramanian\\
Institute of Information and Practical Mathematics\\
Christian-Albrechts-University\\
Kiel, Germany\\
email: ssi@informatik.uni-kiel.de
}
\date{\today}

\maketitle

\begin{abstract}
We give an interpretation of the coefficients of the two variable refinement 
$D_{\Sh_n}(q,t)$  of the 
distance enumerator of the Shi hyperplane arrangement $\Sh_n$ in $n$ dimensions.  This two 
variable refinement was defined by Stanley \cite{stan-rota} for the general $r$-extended Shi
hyperplane arrangements.  We give an interpretation when $r=1$.

We define three natural three-dimensional partitions of the number $(n+1)^{n-1}$.  The
first arises from parking functions of length $n$, the second from special posets (we call
them {\sl tree-posets}) on $n$ vertices defined by Athanasiadis \cite{athan} and the third 
from spanning trees on $n+1$ vertices.  We call the three partitions as the {\sl parking partition},
the {\sl tree-poset partition} and the {\sl spanning-tree partition} respectively.  We show that 
one of the parts of the {\sl parking partition} is identical to the number of edge-labelled trees
with label set $\{1,2,\ldots,n\}$ on $n+1$ unlabelled vertices.  We prove that the 
parking partition majorises the tree-poset partition and conjecture that 
the spanning-tree partition also majorises the tree-poset partition. 

\end{abstract}

\section{Introduction}
Let $r \geq 1$ and $n \geq 2$.  The $r$-extended Shi hyperplane arrangement in 
$n$ dimensions is denoted $\Sh_n^r$.  It is given by the following hyperplanes in $\reals^n$.  

$$
x_i - x_j = -r+1,-r+2,\ldots r, \mbox{for $1\leq i< j \leq n$}
$$

When $r=1$, the arrangement is called the Shi hyperplane arrangement in $n$ dimensions and
denoted $\Sh_n$.  %We will restrict attention to $\Sh_n$ in this work. 
Stanley \cite{stan-rota} defined a two variable distance enumerator of the Shi 
hyperplane arrangement with respect to a base region $B$.  Let $\RSn$ be the set
of regions of the Shi hyperplane arrangement.  Each region $R \in \RSn$ is 
separated from $B$ by a set of hyperplanes and 
let $a$ be the number of separating hyperplanes of the form $x_i - x_j = 0$ 
and $b$ the number of separating
hyperplanes of the form $x_i - x_j = 1$.  The two variable distance enumerator is defined
as $D_{\Sh_n}(q,t) = \sum_{R \in \RSn} q^a t^b$.  
We denote the coefficient of $q^{\ell} t^k$ %as $(k,\ell)$-th entry of the table for 
of $D_{\Sh_n}(q,t)$ as $\Dist_n(k,\ell)$.

\begin{figure}[h]
\centerline{\psfig{file=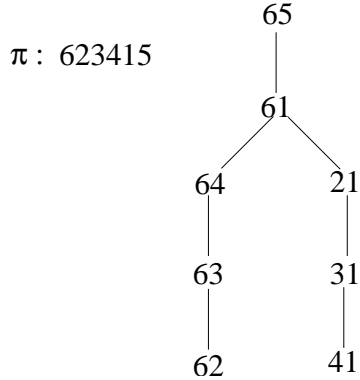,height=5cm}}
\caption{An example of the poset $\IP_{\pi}$}
\label{fig:poset_eg}
\end{figure}

Fix $n$ and let $\Pi_k$ be the set of permutations on $[n]$ which have exactly $k$
non-inversions.  For a permutation $\pi \in \Pi_k$, let $\IP_{\pi}$ be a poset of its
inversions ordered by containment, (ie if $g = (\pi_i,\pi_j)$ where $i < j$, and 
$h = (\pi_a, \pi_b)$ where $a < b$, are inversions, then $g \leq_{\IP_{\pi}} h$ iff
$a \leq i < j \leq b$.  For example, when $\pi = 623415$, the poset $\IP_{623415}$ is 
shown in Figure \ref{fig:poset_eg}.
For $\pi \in \Pi_k$, let the number of {\sl ideals} of $\IP_{\pi}$ with 
${n \choose 2} - k - \ell$ elements be $\IP_{\pi}(\ell)$. 

\begin{theorem}
\label{thm:interpret}
$\Dist_n(k,\ell) = \sum_{\pi \in \Pi_k} \IP_{\pi}(\ell)$
\end{theorem}

Theorem \ref{thm:interpret} gives a two variable generalisation to the equality 
(see Page 96, \cite{EC2}) 
$$\sum_{\pi \in S_n} F(J(\NIP_{\pi}), q) = I_{n+1}(q)$$ 
where $S_n$ is the set of permutations on $n$ distinct alphabets, 
$F(J(\NIP_{\pi}),q)$ is the rank generating function of the lattice of order ideals
of the poset of non-inversion $\NIP_{\pi}$ which is similar to $\IP_{\pi}$, the only 
difference being that we order {\sl non-inversions} of $\pi$ instead of its {\sl inversions}
(please see Remark \ref{rem:embr}).
$I_{n+1}(q)$ is the inversion enumerator of spanning trees on a vertex set of size $n+1$.

\subsection{Three 3d partitions of the regions of $\Sh_n$}
\label{subsec:3d-refin}

For a positive integer $n$, let $[n] = \{1,2,\ldots,n \}$ and let $[n_0] = \{0 \} \cup [n]$.
Let $T$ be a spanning tree on the set $[n_0]$.  We call the vertex 0 as the ``root" of $T$ and
call such trees 0-rooted spanning trees. 

From the bijection between $\RSn$ and spanning trees on $(n+1)$ vertices 
$[n_0] = \{0,1,\ldots,n \}$ (see \cite{stan-lec}), we can view the regions alternatively 
as 0-rooted spanning 
trees on $[n_0]$.  Likewise, we can also view the regions as indexed by Parking
Functions of length $n$. We recall the definition of an $n$ length {\sl parking function}.
There are $n$ parking spaces $0,1,\ldots,n-1$ in a one-way street.  $n$ cars $C_1, C_2,\ldots C_n$ 
enter the street in that order.  $C_i$ has a preferred space $a_i$ and proceeds directly to slot
$a_i$.  If slot $a_i$ is occupied, it will try to park in the next available space.
If a car leaves the street without parking then the process fails.  $\oa=(a_1,a_2,\ldots,a_n)$
is an $n$-length parking function if all cars can park with $a_i$ being their respective choices.
The set of all parking functions of length $n$ is denoted $\PF_n$.  It is known that 
$\oa = (a_1,a_2,\ldots,a_n)$ is a parking function iff the weakly increasing permutation 
$\ob = (b_1,b_2, \ldots, b_n)$ of $\oa$ satisfies $b_i < i$ (see \cite{EC2}).

\subsubsection{Parking Partition}
Let $\oa = (a_1,a_2,\ldots, a_n) \in 
\PF_n$.  
It is simple to check that any permutation of $\oa$ is yields a valid parking function.  
We partition 
$\PF_n$ into the following three parts: those with 
$a_1 > a_2$, with $a_1 = a_2$ and with $a_1 < a_2$.  We call the number of such $n$-length 
parking functions as $gt_n$, $eq_n$ and $lt_n$ respectively.  It is clear that we could have 
chosen 
any indices $i \not= j$ and partitioned $\PF_n$ into three parts as above depending on the 
relation between $a_i$ and $a_j$ and still obtained the same numbers.  Below we tabulate
the numbers $gt_n$, $eq_n$ and $lt_n$ for small values of $n$.

\subsubsection{Tree-poset partitions}
We define the {\sl tree-poset partition} next. 
Consider the hyperplane $x_1 - x_2 = \alpha$ for 
$\alpha = 0, 1$; and let $R \in \RSn$.  
Let $\oa_R = (a_1, a_2, \ldots, a_n)$ be any point in $R$.  Clearly, the value 
$a_1 - a_2$ is either $<0$, strictly between 0 and 1, or $>1$ and this 
condition is independent of the point $\oa_R$.  Thus each region $R$ with respect to the
dimensions $x_1$ and $x_2$ satisfies one of the three properties: all points $\oa_R \in R$
either have $a_1 - a_2 < 0$, or $0 < a_1 -a_2 < 1$ or $a_1 - a_2 > 1$.

Let $R_n^{\> \cdot < 0}$, $R_n^{0 < \cdot < 1}$  and $R_n^{\> \cdot > 1}$ respectively
denote the number of regions satisfying the above three conditions.
The main reason for this definition is to
understand how $\RSn$ gets partitioned by the parallel hyperplanes $x_1 - x_2 = 0,1$.

%It is clear that $R_n^{\> \cdot < 0} + R_n^{0 < \cdot < 1} +  R_n^{\> \cdot > 1} = (n+1)^{n-1}$.  
Below we tabulate the numbers  $R_n^{\> \cdot < 0}$, $R_n^{0 < \cdot < 1}$ and 
$R_n^{\> \cdot > 1}$.  For this definition, the numbers are not necessarily 
independent of the choices 1 and 2.

\subsubsection{Spanning-tree partitions}
Lastly, we define the {\sl spanning-tree partition}.  Let $v_1,v_2 \in [n]$, $v_1 \not= v_2$ 
be two fixed vertices, and let $T$ be 
a 0-rooted spanning tree on $[n_0]$.  There are again three possibilities for the following 
path relation:  
either $v1$ is on the unique $v2$-0 path; or $v2$ is in the unique $v1$-0 path; or neither
of the two happens.   Let $T_n^{v1}$, $T_n^{v2}$ and $T_n^{disj}$ be the number of 0-rooted
spanning trees on $[n_0]$ for each of the above three choices.  These numbers are again
independent of the choices $v_1,v_2$.
We tabulate the numbers $T_n^{disj}$, $T_n^{v1}$ and $T_n^{v2}$ for small
values of $n$ below.  

\begin{center}

\begin{tabular}{ccccc} 
\begin{tabular}{|r|r|r|r|} \hline
$n$ & $gt_n$ & $lt_n$ & $eq_n$ \\ \hline \hline
3 & 6 & 6 & 4 \\ \hline
4 & 50 & 50 & 25 \\ \hline
5 & 540 & 540 & 216  \\ \hline
6 & 7203 & 7203 & 2401  \\ \hline
\end{tabular}

& \hspace{.2 cm} &
\begin{tabular}{|r|r|r|r|} \hline
$n$ & $R_n^{\> \cdot < 0}$ & $R_n^{0 < \cdot < 1}$ & $R_n^{\> \cdot > 1}$ \\ \hline \hline

3 & 6 & 5 & 5 \\ \hline
4 & 50 & 37 & 38 \\ \hline
5 & 540 & 366 & 390  \\ \hline
6 & 7203 & 4553 & 5051 \\ \hline
\end{tabular} 

& \hspace{.2 cm}  & 

\begin{tabular}{|r|r|r|r|} \hline
$n_0$ & $T_n^{disj}$ & $T_n^{v1}$ & $T_n^{v2}$ \\ \hline \hline
3 & 6 & 5 & 5 \\ \hline
4 & 51 & 37 & 37 \\ \hline
5 & 564 & 366 & 366  \\ \hline
6 & 7701 & 4553 & 4553  \\ \hline
\end{tabular} 
\end{tabular}

\end{center}

For $n \geq 1$, let 
$UT_n$ be the number of edge labelled trees with label set $\{1,2,\ldots,n \}$ on 
$n+1$ unlabelled vertices.  It is known (see \cite{EC2}) that $UT_n = (n+1)^{n-2}$.  

\begin{theorem}
\label{thm:equality}
For all $n \geq 1$, $eq_n = UT_n$.
\end{theorem}

We show the following majorisation theorem.

\begin{theorem}
\label{thm:park-major-poset}
For $n \geq 2$, the largest part of the {\sl parking partition} is equal to the 
largest part of the {\sl tree-poset partition}. Hence,
the sorted {\sl parking partition} majorises the sorted 
{\sl tree-poset partition}.  
\end{theorem}

\section{Two variable distance enumerator: an interpretation}
In this section, we prove Theorem \ref{thm:interpret}.  We use a
poset representation for each region $R \in \RSn$.  This representation was defined
by Athanasiadis \cite{athan}.

\subsection{The posets of Athanasiadis}

\begin{figure}[h]
\centerline{\psfig{file=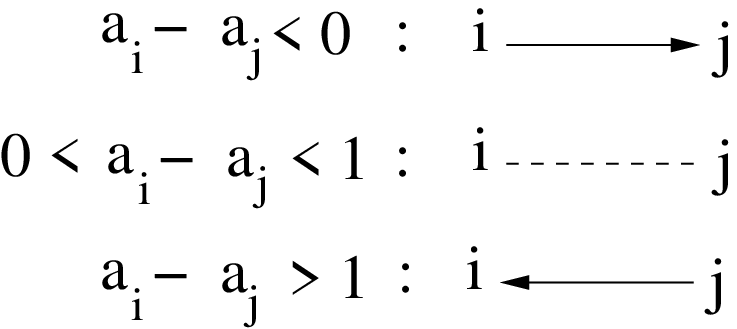,height=2.5cm}}
\caption{Representing the three possibilities, where $i < j$.}
\label{fig:poset}
\end{figure}

Let $\oa_R$ be a point of $R \in \RSn$.  Represent each of the 
three possibilities $a_1 - a_2 < 0$, $0 < a_1 - a_2 < 1$ and $a_1 - a_2 > 1$  
Figure \ref{fig:poset} (the dotted lines in the second figure 
represent an incomparability relation between the vertices $i$ and $j$). 

We call arcs of the form $(i,j)$ where $i<j$ as {\sl forward arcs} and those of the form
$(j,i)$ where $i<j$ as {\sl backward arcs}. 

\begin{figure}[h]
\centerline{\psfig{file=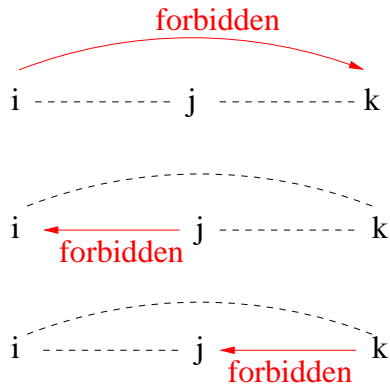,height=5cm}}
\caption{The three forbidden subposets where $i < j < k$.}
\label{fig:forbid}
\end{figure}

Athanasiadis 
\cite{athan} showed that this representation yields a poset on $[n]$ and that such
posets do not have the three subposets shown in Figure \ref{fig:forbid}.  
%This gives a poset representation of each region of $\RSn$.  
Athanasiadis also proved that
any poset without these three ``forbidden" subposets arose from a region thereby
characterising such posets. 
We refer to such posets as ``tree-posets".

\begin{lemma}
\label{lem:posets}
$\Dist_n(k,\ell)$ is equal to the number of tree posets on $[n]$ which
have $k$ forward arcs and $\ell$ backward arcs.   
\end{lemma}

\begin{proof}
%Let $n$ be the dimension of hyperplane arrangement.  
We use the Pak and Stanley method of starting from a base region $B$, assigning a poset
$P_B$ to it and use their rules to get a poset for adjacent regions.  We recall that
the base region is the region bounded by the hyperplanes 
$x_1 > x_2 > \cdots > x_n > x_1 - 1$.  Assign the $n$ element antichain to this region.
When region $R'$ is separated from $R$ by the hyperplane $x_i - x_j = 0$ (for $i<j$), 
set $P_{R'} = P_{R} \cup (i,j)$, ie add the {\sl forward} arc $(i,j)$.  Similarly, when 
region $R'$ is separated from $R$ by the hyperplane $x_i - x_j = 1$ (and $i< j$), then
set $P_{R'} = P_{R} \cup (j,i)$, ie add the {\sl backward} arc $(j,i)$.  We note that 
whenever we cross from region $R$ %to $R'$ 
in this manner, we always cross a 
hyperplane $x_i - x_j = 0,1$ such that, in $P_R$ the vertices $i$ and $j$ are 
incomparable. Thus, the above two cases are exhaustive.  The proof of Pak and Stanley 
shows that this algorithm is well defined over different shortest paths from the base
region $B$ to any other region $R$.

This is the same labelling given by Athanasiadis \cite{athan}, though the algorithm 
explains the proof better.  From the above algorithm, we see the following invariant:
when we cross the hyperplane $x_i - x_j = 0$ (for $i < j$), we get a forward arc 
$(i,j)$ and when we 
cross the hyperplane $x_i - x_j = 1$, we get a backward arc $(j,i)$.  Thus a region
$R$ contributes to $\Dist_n(k,\ell)$ iff $P_{R}$ has $k$ forward
arcs and $\ell$ backward arcs.
\end{proof}

We fix the dimension $n$ and call such tree-posets (on $[n]$) with $k$ 
forward arcs and $\ell$ backward arcs as $(k,\ell)$-tree posets.  We note that by 
this propagation, we eventually obtain posets $P_R$ that are permutations.  It is 
simple to see that {\sl backward arcs} in permutations correspond to its {\sl 
inversions}.  We use both {\sl inversions} and {\sl backward arcs} interchangeably 
even when the poset is not a permutation (ie has incomparable elements).
It is also simple to observe that such posets
$P_R$ with no {\sl inversions} (ie only forward arcs and incomparability relations) 
are the ``nearest" regions to the base region $B$ in the regions of $\B_n$ (the Braid
hyperplane arrangement).

A reverse propagation shown below proves that we can start from $(k,\ell+1)$-tree posets 
and by converting an inversion into an
incomparability relation, obtain all $(k,\ell)$-tree posets.

\begin{lemma}
All $(k,\ell)$-tree posets can be obtained from $(k,\ell+1)$-tree 
posets by converting a backward arc into an incomparability relation.  
\end{lemma}

\begin{figure}[h]
\centerline{\psfig{file=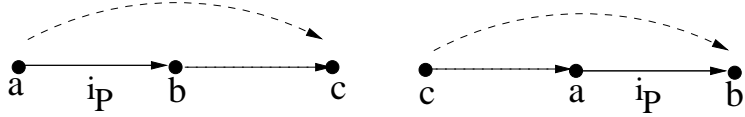,height=1.5cm}}
\caption{Contradicting the maximality of $i_P$ .}
\label{fig:contrad}
\end{figure}

\begin{proof}
Let $P$ be a $(k,\ell)$-tree poset on $[n]$.  We prove this by induction on 
${n \choose 2} - (k + \ell)$.  The base case when $k+ \ell = {n \choose 2}$ 
ie when $P$ is a permutation, is simple.  
%The statement is vacuously true as there are no $(k,\ell)$-tree posets with $k+ \ell > {n \choose 2}$.  

Let $P$ be a $(k,\ell)$-tree poset with $k + \ell < {n \choose 2}$.  We exhibit a
$(k, \ell+1)$-tree poset $Q$ and identify an inversion $i_P$ in $Q$ such that
$P = Q - \{i_P \}$.  (ie We convert one incomparability relation $i_P$, (we
also call these as non-arcs)  in $P$ into a backward arc and get a poset $Q$ which
does not have any of the three forbidden subposets.)  Clearly $P = Q - \{i_P\}$ and 
by induction, we will be done.

\begin{figure}[h]
\centerline{\psfig{file=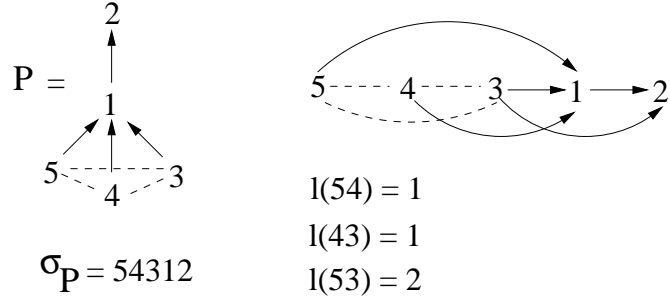,height=4cm}}
\caption{Linear extensions with the ``largest" vertex first and lengths $l(g)$.}
\label{fig:linextn}
\end{figure}

Let $P$ be a $(k,\ell)$-tree poset. Consider the linear extension $\sigma_P$ of $P$ 
where we 
break ties when they exist, by the ``largest" vertex first rule (see Figure 
\ref{fig:linextn} for an 
example).  Order the vertices (ie $[n]$) according to $\sigma_P$ so that all 
arcs of $P$ are directed towards the right and 
let $g = (\sigma_r, \sigma_s)$ for $r < s$ be a non-arc in $P$.
Define the ``length" $l(g)$ of $g$ to be the number of non-arcs contained within it 
(including itself) in $\sigma_P$ (ie the number of non-arcs in the subpermutation 
$\sigma[r,s] = (\sigma_r, \sigma_{r+1},\ldots, \sigma_s)$.  Let $i_P$ be a non-arc of maximum 
length.  We claim that we can convert $i_P$ into a backward arc and obtain a 
$(k,\ell+1)$-tree poset, $Q$.

To prove this, we need to show that none of the three forbidden subposets 
appear in $Q$.  Since they do not appear in $P$, $i_P$ must be involved in
any forbidden subposet.  Suppose the second or the third forbidden subposet of Figure 
\ref{fig:forbid} appears in $Q$.  Then, $i_P$ must be the backward-arc in a forbidden
subposet and we have a situation shown in Figure \ref{fig:contrad}.  
In both cases, let $i_P = (a,b)$.  We note that 
$a > b$.  In the first case, there is a non arc $(b,c)$ and since $b$ appears to the
left of $c$ in $\sigma_P$, $b>c$.  Thus we have $a> b > c$ and this violates the maximality
of $i_P = (a,b)$.  In the second case, similarly there are three vertices $c > a > b$ 
and the maximality of $i_P$ is again violated.
The first forbidden subposet cannot appear in $Q$ as we need a forward arc for it and
we convert $i_P$ into a backward arc.
Thus $Q$ has one less non-arc, has no forbidden subposets 
and (hence by the Theorem 1.1 \cite{athan}) corresponds to a region of $\Sh_n$. 
This completes the proof.
\end{proof}

%\end{proof}

Let $n$ be fixed and for $0 \leq k \leq {n \choose 2}$, let $\Pi_k$ be the set of
permutations $\pi$ on $[n]$ having ${n \choose 2} - k$ inversions.  For 
$\pi \in \Pi_k$ and $0 \leq \ell \leq {n \choose 2} - k$, let the number of 
$(k,\ell)$-tree posets obtained from $\pi$ (ie those obtained from $\pi$ by
deleting ${n \choose 2} - k - \ell$ inversions) be denoted $\pi(k,\ell)$.

\begin{corollary}
With the above notation, 
$\Dist_n(k,\ell) = \sum_{\pi \in \Pi_k} \pi(k,\ell)$
\end{corollary}
\begin{proof}
It is simple to check that $(k,\ell)$-tree posets obtained from $\sigma 
\not= \pi, \pi, \sigma \in \Pi_k$ are different and that we can just add up the 
numbers $\pi(k,\ell)$ over different $\pi \in \Pi_k$.
\end{proof}
Thus to get $(k,\ell)$-tree posets, we could start from $\pi \in \Pi_k$ and delete
${n \choose 2} - k - \ell$ inversions such that the three forbidden posets do not occur.

\begin{example}
We show an example of the inversion-deletion process described above.  
Let $\pi = 623415$.  We can 
order the $(k,\ell)$-tree posets we obtain by containment.  This poset of
$(7,\ell)$-tree posets obtained from $\pi$ is shown in Figure \ref{fig:lattice_eg}.  The 
edges of the poset are all oriented rightwards and the edge labels are the inversions 
converted into non-arcs from the (poset corresponding to the) previous vertex.

\begin{figure}[h]
\centerline{\psfig{file=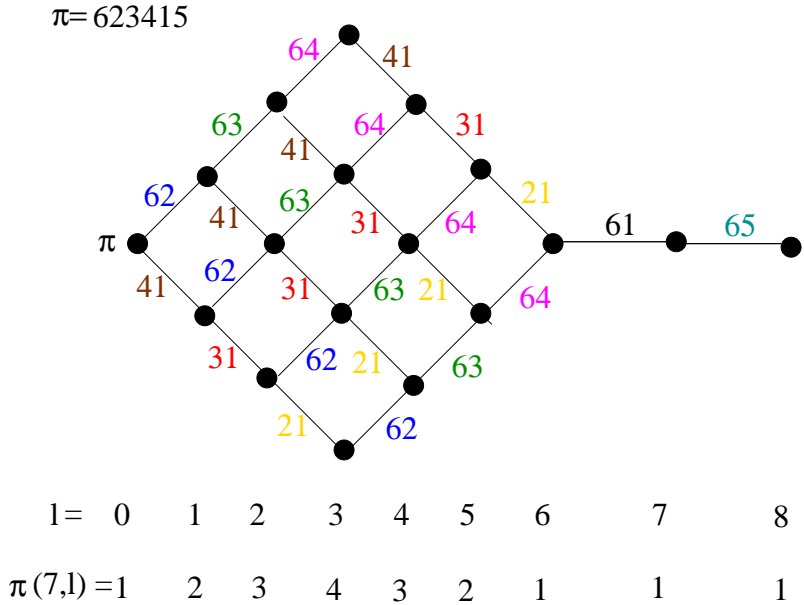,height=8cm}}
\caption{Example of the poset of $(7,\ell)$-posets for $0 \leq \ell \leq 8$  }
\label{fig:lattice_eg}
\end{figure}
\end{example}

\begin{remark}
All $(k,\ell)$-tree posets arising from $\pi \in \Pi_k$ have $\pi$ as a linear
extension and when we break ties due to incomparability using the ``largest vertex"
first rule, these posets have $\pi$ as the linear extension.  Further, all points
$\alpha_R = (\alpha_1, \alpha_2, \ldots, \alpha_n)$ in a region $R$ corresponding 
to a $(k,\ell)$-tree poset from $\pi \in \Pi_k$  have
$\pi$ as the permutation when the $\alpha_i$'s are sorted in increasing order.

\end{remark}

\subsection{Ideals of the inversion poset of a permutation}
We give an alternate interpretation of the number of $(k,\ell)$-tree posets on $[n]$.  
It is clear that inversions of permutations $\pi \in \Pi_k$ are to be deleted in 
some sequence to obtain $(k,\ell)$-tree posets.  Such sequences are described below.

Let $\pi =(\pi_1,\pi_2\ldots,\pi_n) \in \Pi_k$.  
Let the sub-permutation between two indices $i < j$ be denoted $\pi[i,j]$, 
ie $\pi[i,j] = (\pi_i, \pi_{i+1}, \ldots, \pi_j)$.
Let $p = (\pi_i,\pi_j)$ be an inversion.  Let $\inn_p$ be the number of inversions 
in $\pi[i,j]$.

\begin{lemma}
\label{lem:inner}
Let $\pi \in \Pi_k$ and let $g$ be an inversion in $\pi$.  In the deletion process
described above, $g$ can be converted into a non-arc only after all the inversions 
strictly within it have been converted.
\end{lemma}

\begin{proof}
Let $P$ be a poset obtained from $\pi$ and suppose we could convert an inversion 
$(a,b)$ to a non-arc while an inner inversion $(c,d)$ remained (ie if 
$\pi_w = a, \pi_x = b, \pi_y = c$ and $\pi_z = d$, we have $w < y < z <x$, for 
example, see Fig \ref{fig:inner}).  

\begin{figure}[h]
\centerline{\psfig{file=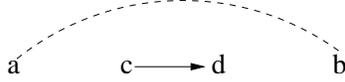,height=1cm}}
\caption{Inner inversions must be deleted earlier}
\label{fig:inner}
\end{figure}

If both $(a,c)$ and $(d,b)$ were arcs, then the posetness of $Q$ would
be violated.  Since all arcs go rightwards, we assume $(a,c)$ is an incomparable
pair ie that $a > c$ ie this was an inversion that got converted.  Since $(c,d)$ is
an inversion, $c>d$.  Thus $a >d$ and this inversion either stays as an inversion
or has been converted into a non-arc.  If it is a non-arc, we have a forbidden 
subposet on $a,c,d$ and hence $(a,d)$ remains as an arc.  If $(d,b)$ exists, then
again we violate posetness and hence $(d,b)$ is a non-arc and this induces a
forbidden subposet on $a,d,b$.  The argument is identical if we had started
with $(d,b)$ being a non-arc.
\end{proof}

\begin{proof} (Of Theorem \ref{thm:interpret})
From Lemma \ref{lem:inner}, we see that for $\pi \in \Pi_k$, the ideals of 
$\IP_{\pi}$ with ${n \choose 2} -k - \ell$
elements are precisely the elements constituting $\pi(k,\ell)$.  Lemma
\ref{lem:posets} completes the proof.
\end{proof}

We note that the earlier poset obtained by inversion-deletion is actually a distributive 
lattice and that it is isomorphic to the lattice of order ideals $J(\IP_{\pi})$ (where
$\pi$ is the starting permutation).

\begin{remark}
\label{rem:embr}
We are essentially assigning two values to each ``embroidered permutation" 
(see Page 81, \cite{stan-lec}), though we use inversions instead of 
non-inversions. The region of $\RSn$ that $(\pi, {\cal C})$ represents is slightly
different for us.  Suppose $\pi  = (\pi_1,\pi_2,\ldots,\pi_n)$, then the region 
corresponding to this embroidered permutation is 
$x_{\pi_1} < x_{\pi_2} < \cdots < x_{\pi_n}$ and $\forall g=(i,j) \in {\cal C}$,
$0 < x_j - x_i < 1$.  This is why we need $g$ to be an inversion.  
We are assigning two parameters 
$(a,b)$ to each embroidered permutation $(\pi,{\cal C})$ where $a$ is the number of 
non-inversions of $\pi$ and $b$ is the total number of inversions contained in the
family ${\cal C}$.  
\end{remark}

\begin{remark}
Let $\pi \in \Pi_k$.  Because there is a single hyperplane separating regions corresponding
to $(k,\ell)$-tree posets and $(k,\ell-1)$-tree posets, the lattice $J(\IP_{\pi})$ when 
treated as a graph is the subgraph of
distance graph of $\RSn$ with respect to the base region $B$ consisting of those 
regions of $\RSn$ which sit inside a given region of $\B_n$ (the Braid arrangement).
\end{remark}

\section{Results on the 3d partitions}

We collect some properties of each of the partitions below.

\subsection{Properties of the partitions}
\label{subsec:order}

We prove some properties about the order of the components of  
the three 3d partitions.

\begin{lemma}
For $n \geq 2$, the {\sl parking partition} satisfies $gt_n = lt_n \geq eq_n$.
\end{lemma}
\begin{proof}
It is known that $\oa = (a_1,a_2,\ldots,a_n)$ is a parking function iff its weakly increasing
permutation $\ob = (b_1,b_2,\ldots,b_n)$ satisfies the relation $b_i < i$.  Let 
$\oa = (a_1,a_2,\ldots,a_n) \in \PF_n$ with $a_1 > a_2$.  Clearly, $\oa' = (a_2,a_1,\ldots,a_n)$ 
obtained from $\oa$ by swapping the first two coordinates is also a valid parking function, and
has $a_1' < a_2'$.  The argument is reversible and this bijection proves that $gt_n = lt_n$.

We show that $lt_n \geq eq_n$.  Let $\oa \in eq_n$.  
Let $\ob = (a_1, a_2+1, a_3,\ldots,a_n)$ and $\oc = (c_1, c_2, \cdots, c_n)$ be a weakly 
increasing permutation of $\ob$. 
We show that $\ob \in lt_n$.  We only need to check that $\ob \in \PF_n$.  Suppose not, then 
there is an index $k$ such that $c_k \not< k$.  Since we changed only one coordinate to obtain
$\ob$ from $\oa$, $c_k = a_2+1$.  But then $a_1 = a_2-1$ will be  $c_{k-x}$ for $x \geq 1$ 
and thus we get $\oa \not \in \PF_n$ which is a contradiction.

\end{proof}

\begin{lemma}
\label{lem:poset-partn}
For $n \geq 2$, the {\sl tree-poset partition} satisfies 
$R_n^{\> \cdot < 0} \geq  \max( R_n^{\> \cdot > 1} , R_n^{0 < \cdot < 1})$. 
\end{lemma}
\begin{proof}
We first prove that $R_n^{\> \cdot < 0} \geq R_n^{\> \cdot > 1}$.  To do this, we note that
by Theorem \ref{thm:interpret}, the regions $R_n^{\> \cdot < 0}$ are those which 
have $(1,2)$ as a forward arc and the regions
of $R_n^{\> \cdot > 1}$ are those which have $(1,2)$ as a backward arc, with the condition
that $(1,2)$ has not been converted into an incomparability relation.

We will show a slightly stronger property:  consider all permutations $\pi$ of $[n]$ in which
1 precedes 2 (ie $(1,2)$ is a forward arc).  Such permutations contribute $|J(\IP_{\pi})|$
elements to $R_n^{\> \cdot < 0}$ and only such permutations contribute
to $R_n^{\> \cdot < 0}$.  

For each such $\pi$, let $\pi'$ be the permutation obtained by inverting the position of the 
elements 1 and 2.  Similar to the above argument, every region of $R_n^{\> \cdot > 1}$ 
occurs from $\pi'$ and an
ideal of $\IP_{\pi'}$ which does not contain the inversion $\{2,1\}$ (and hence all elements
$X = \{x \geq_{\IP_{\pi'}} \{2,1\} \}$).  Let $\IP_{\pi'}(21)$ denote the subposet 
$\IP_{\pi'} - X$.  It is simple to see that the $\IP_{\pi'}(21)$ is a subposet of $\IP_{\pi}$ 
as well.  Thus the number of order ideals is smaller for each $(\pi', \pi)$ pair and summing
over these pairs completes the proof.

An almost identical proof works to shows that $R_n^{\> \cdot < 0} \geq R_n^{0 < \cdot < 1}$.
We note that $R_n^{0 < \cdot < 1}$ is the number of $(\pi',{\cal I})$ pairs where $\pi'$ is 
a permutation with $2$ preceding $1$ and ${\cal I}$ is an ideal of $\IP_{\pi'}$ such that
the inversion $(2,1) \in {\cal I}$. Thus $X = \{x| x <_{\IP_{\pi'}} (2,1) \} \in {\cal I}$
as well.  Let $\IP_{\pi'}(2,1) = \IP_{\pi'} - X$.  The remaining argument is identical.  
\end{proof}

\begin{lemma}
\label{lem:sptree-part-order}
For $n \geq 2$, the {\sl spanning-tree partition} satisfies 
$T_n^{disj} \geq T_n^{v1} = T_n^{v2}$ .
\end{lemma}
\begin{proof}
We first prove that $T_n^{v1} = T_n^{v2}$.  Let $T \in T_n^{v1}$. Thus $T$ is a 0-rooted spanning
tree on $[n_0]$ and $v1$ is on the unique $v2-0$ path.  By swapping the vertices $v2$ and $v1$, we
get a tree $T' \in T_n^{v2}$.  The equality part of the Lemma is thus proved.  
 
To show that $T_n^{disj} \geq T_n^{v1}$, let $T \in T_n^{v1}$ as before. Let $T''$ be obtained from
$T$ by swapping $v1$ and 0.  Clearly $T'' \in T_n^{disj}$.  %This completes the proof of the Lemma.  
\end{proof}

\subsection{Properties among the partitions}
In this section, we prove Theorems \ref{thm:equality} and \ref{thm:park-major-poset}.
We recall that $eq_n$ is the number of $\oa \in \PF_n$ which satisfy $a_1 = a_2$.  

\begin{proof}(Of Theorem \ref{thm:equality})
The proof of Pollack given in \cite{stan-lec}(Page 92) to count the number of $n$-length 
Parking functions carries
over exactly.  
\end{proof}

\begin{proof}(Of Theorem \ref{thm:park-major-poset})
We first show that the largest elements the {\sl parking 
partition} is equal to the largest element of the {\sl tree-poset partition}.  Since
the partitions are 3-dimensional, this is sufficient to prove the majorisation result.

We use the bijection of Pak and Stanley \cite{stan-rota}, coupled with the forbidden 
subposets of Athanasiadis \cite{athan}.  By Lemma \ref{lem:poset-partn}, 
$R_n^{\> \cdot < 0}$ is the largest part of the 
tree-poset partition.  We recall that the posets $P_R$ of such a region $R$ has a
forward arc $(1,2)$ between vertices 1 and 2. 

We first show that when the poset $P_R$ of a region $R$ has $(1,2)$ as a forward
arc, then the corresponding 
parking function $\oa_R$ of $R$ under the bijection of Pak and Stanley has the property 
$a_1 > a_2$.  
Since $(1,2)$ is a forward arc, 
$R$ is on the ``less than" side of the hyperplane $x_1 - x_2 = 0$.  Since the
{\sl base region}  $B$ has $x_1 > x_2$, we must cross the hyperplane $x_1 - x_2 = 0$ at some 
point in any shortest
distance path from $B$ to $R$.  This crossover will contribute a 1 to $a_1$, the first component of the 
parking function $\oa$ and 0 to $a_2$.  It is simple to check that the only way to increase $a_2$ is 
to cross the 
hyperplane $x_2 - x_v = 0$ for some $v \in [n] - \{1,2\}$ on a path from $B$ to $R$.  
All such crossovers are recorded by a forward arc $(2,v)$ in the poset representation of
$R$.  For such vertices $v$, since $(2, v)$ and $(1, 2)$ are forward arcs, by transitivity of the
poset, $(1, v)$ is also a forward arc and this means we contribute a 1 to $a_1$ as well. 
This completes the proof of one half of the bijection.

For the other half, let $\oa \in gt_n$.  We claim that its corresponding region $R$ under the 
bijection of Pak and Stanley has $(1,2)$ as a forward arc.  As before, if $a_2 = k$, there 
exists a set
$S$ with $|S| = k$ such that for all $v \in S$, $(2,v)$ is a forward arc.  Similarly, 
when $a_1 = k+x$ for
$x > 0$, there is a set $T$ such that for all $v \in T$, $(1,v)$ is a forward arc.  We claim 
that $2 \in T$.
Suppose not, then there is a vertex $v \in T - S$, $v \not = 2$ such that $(1,v)$ is a forward
arc and 
$(2, v)$ is not (see Figure \ref{fig:region}).  Thus there are two cases for the relation 
between 2 and $v$.

\begin{figure}[h]
\centerline{\psfig{file=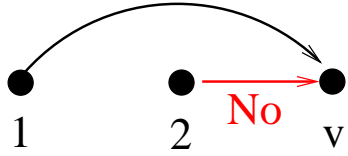,height=2cm}}
\caption{When $(1,v)$ is a forward arc and $(2,v)$ is not.}
\label{fig:region}
\end{figure}

\begin{itemize}
\item When $(v,2) $ is a forward arc :   As $(1,v)$ and $(v,2)$ are forward arcs, by transitivity 
$(1,2)$ too is.

\item When $(2,v)$ is an incomparability : If $(1,2)$ is a backward arc, then transitivity 
among these
three vertices would be violated.  If $(1,2)$ were an incomparability relation, then we would get
the first forbidden subposet of Figure \ref{fig:forbid} on the vertices $1,2,v$.

\end{itemize}
This completes the proof of the theorem.
\end{proof}

\begin{conjecture}
Similar to Theorem \ref{thm:park-major-poset}, the smallest parts of spanning-tree partition
and the tree-poset partition are equal.
For $n \geq 2$, the sorted {\sl spanning-tree partition} majorises the sorted 
{\sl tree-poset partition}.
\end{conjecture}

\begin{conjecture}
For fixed $n,k$, the numbers $\Dist_n(k,\ell)$ as $\ell$ increases are unimodal.  
\end{conjecture}

\begin{question}
Is there a recurrence or a generating function for the numbers occuring in 
the spanning tree partition?
\end{question}

\bibliography{enumerate}
\end{document}